\documentclass[11pt]{amsart}

\usepackage{amsmath,amsfonts,amscd,amssymb,epsfig,euscript, eufrak}
\usepackage{amsxtra,mathrsfs}
\newtheorem{theorem}{Theorem}
\newtheorem{proposition}[theorem]{Proposition}
\newtheorem{lemma}[theorem]{Lemma}

\theoremstyle{definition}
\newtheorem{definition}[theorem]{Definition}

\theoremstyle{remark} \newtheorem{remark}[theorem]{Remark}
\numberwithin{equation}{section}
\newcommand{\field}[1]{\ensuremath{\mathbb{#1}}}
\newcommand{\CC}{\field{C}}

 \DeclareMathOperator{\re}{Re}

\newcommand{\del}{\partial}
\newcommand{\delb}{\bar\partial}

\newcommand{\vep}{\varepsilon}

\newcommand{\z}{\bar{z}}
\newcommand{\w}{\bar{w}}

\begin{document}
\title[Orthogonal polynomials and $\delb$-problem]{Normal matrix models, $\delb$-problem, and orthogonal polynomials on the complex plane}
\author{Alexander R. Its }
\address{Department of Mathematical Sciences \\
Indiana University-Purdue University Indianapolis\\  
Indianapolis, IN 46202-3216 \\ USA}
\email{itsa@math.iupui.edu}

\author{Leon A. Takhtajan}
\address{Department of Mathematics \\
Stony Brook University\\ Stony Brook, NY 11794-3651 \\ USA}
\email{leontak@math.sunysb.edu}
\maketitle

\begin{abstract} We introduce a  $\delb$-formulation of
the  orthogonal polynomials on the complex plane, and hence
of the related normal matrix model,  which is expected to play the same role as the Riemann-Hilbert formalism in the theory of orthogonal polynomials on the line and for the related Hermitian model. We propose an analog of Deift-Kriecherbauer-McLaughlin-Venakides-Zhou asymptotic method for the analysis of the relevant
$\delb$-problem, and indicate how familiar steps for the Hermitian model, e.g. the  $g$-function ``undressing'', might look like in the case of the normal model. We use the particular model considered recently by  P. Elbau and G. Felder as a case study.
\end{abstract}

\section{Introduction}
In these notes we attempt to develop for the normal matrix model
a formalism analogous to the Riemann-Hilbert method in the
theory of Hermitian matrix model. As in the latter case, the starting point is proper analytical  characterization of the 
relevant orthogonal polynomials. Unlike the Hermitian
matrix model, the orthogonality condition for the polynomials
associated with the  normal model is formulated
with respect to a measure on the plane. This, as we will
see below, leads to the replacement of the Riemann-Hilbert
problem of \cite{FIK1} by a certain $\delb$-problem. We shall 
present in detail the setting of the $\delb$-problem
for the case of what we will call in these notes the
{\it Elbau-Felder  model}. This model arises as a natural
regularization of the normal matrix model of
P. Wiegmann and A. Zabrodin in \cite{WZ,KKMWWZ}
by restricting the matrix integral of the latter to  normal matrices
whose  eigenvalues lie in a compact domain $D$ of the complex plane. 
Using the  Elbau-Felder  model as a case study, we shall also outline a possible 
 $\delb$-version of the Deift-Kriecherbauer-McLaughlin-Venakides-Zhou (DKMVZ)  asymptotic model.  The DKMVZ method proved to be very efficient in the asymptotic
 analysis of the oscillatory Riemann-Hilbert problems appearing in 
 Hermitian  matrix model.  We  have not yet
 succeeded in providing complete generalization of the DKMVZ scheme
 for the orthogonal polynomials on the plane;
 in fact, we rather highlighted the challenging difficulties to be
 overcome. We hope, however,  that these notes might stimulate further development of the analog DKMVZ asymptotic method  for the orthogonal
 polynomials on the plane and related normal matrix models.

\section{Preliminaries}
 \subsection{Normal matrix models and orthogonal polynomials} 
  Let $D$ be a bounded domain on the complex plane $\CC$ containing the origin, and let $V(z)$ be a real-valued smooth function on
 $\CC$. Following P. Elbau and G. Felder \cite{EF}, we shall study the 
 normal matrtix model  characterized by the partition function  $Z_{N}$ defined
 by the following $N$-fold integral,
 \begin{equation*}
 Z_{N}=\idotsint\limits_{D^{N}}\prod_{i\neq j}|z_{i}-z_{j}|^{2}e^{-N\sum_{k=1}^{N}V(z_{k})}d^{2}z_{1}\cdots d^{2}z_{N}.
 \end{equation*}
 
  Let $\chi_{D}$ be the characteristic function of the domain $D$. 
 \begin{definition} Orthogonal polynomials on $\CC$ with respect to the measure $e^{-NV(z)}\chi_{D}(z)d^{2}z$ are polynomials $P_{n}(z)=z^{n}+a_{n-1n}z^{n-1}+\cdots+a_{0n}$, satisfying
\begin{equation} \label{orthogonality}
\iint\limits_{D}P_{n}(z)\overline{P_{m}(z)}e^{-NV(z)}d^{2}z=h_{n}\delta_{mn}\quad\text{for all}\quad m,n=0,1, 2\dots.
\end{equation}
 \end{definition}
The following lemma is standard.
 \begin{lemma}\label{lemma1} 
 \begin{equation*}
 Z_{N}=N!\prod_{n=0}^{N-1}h_{n}.
 \end{equation*}
 \end{lemma}
 The proof  is exactly
 the same as in the case of the Hermitian model (see e.g. \cite{D2}). As in the
 case of the Hermitian model, Lemma \ref{lemma1} reduces the question
 of the asymptotic analysis of the partition function
 $Z_{N}$ as $N\rightarrow\infty$ to the asymptotic analysis of the orthogonal
 polynomials $P_{n}(z)$ as $n, N \to \infty$.
\section{Matrix $\delb$-problem}
Using the orthogonal polynomials on the line as an analogy (see \cite{FIK1}, \cite{D2}), we  set 
 \begin{equation} \label{matrix-Y}
 Y_{n}(z)=\begin{pmatrix} 
 P_{n}(z) & \frac{1}{\pi}\iint\limits_{D}\frac{\overline{P_{n}(z')}}{z'-z}e^{-NV(z')}d^{2}z'\\
- \frac{\pi}{h_{n-1}}P_{n-1}(z) &  -\frac{1}{h_{n-1}}\iint\limits_{D}\frac{\overline{P_{n-1}(z')}}{z'-z}e^{-NV(z')}d^{2}z'
 \end{pmatrix}.
 \end{equation}
It follows from the formula
 $$\frac{\del}{\del\z}\,\frac{1}{z-z'}=\pi\delta(z-z'),$$
 understood in the distributional sense, that
 \begin{equation} \label{d-bar-Y}
 \frac{\del }{\del\z}Y_{n}(z)=\overline{Y_{n}(z)}(I-G(z)),
 \end{equation}
 where $I$ is $2\times 2$ identity matrix and
 \begin{equation} \label{matrix-G}
 G(z)=\begin{pmatrix} 1 & e^{-NV(z)}\chi_{D}(z) \\
 0 & 1
 \end{pmatrix}
 \end{equation}
 The following proposition is central (cf. the case of the orthogonal polynomials on the line).
 \begin{proposition}\label{prop1}
  The matrix $Y_{n}(z)$ is the unique solution of the $\delb$-problem \eqref{d-bar-Y}--\eqref{matrix-G} with the normalization
 \begin{equation} \label{asymptotics-Y}
 Y_{n}(z)=\left(I+O\left(\frac{1}{z}\right)\right)\begin{pmatrix} z^{n} & 0 \\
 0 & z^{-n}\end{pmatrix} \equiv  \left(I+O\left(\frac{1}{z}\right)\right)z^{n\sigma_{3}},
 \end{equation}
 as $|z|\rightarrow\infty$,
 where $\sigma_{3}=(\begin{smallmatrix} 1 & 0\\
 0 & -1\end{smallmatrix})$.
 \end{proposition}
 \begin{proof}
 It follows from the geometric series expansion
 $$\frac{1}{z-z'}=\frac{1}{z}\sum_{k=0}^{\infty}\left(\frac{z'}{z}\right)^{k}$$
 as $|z|\rightarrow\infty$, and the property \eqref{orthogonality}, rewritten as
\begin{equation} \label{orthogonality-2}
\iint\limits_{D}P_{n}(z)\z^{m}e^{-NV(z)}d^{2}z=h_{n}\delta_{mn},
\end{equation}
that the matrices \eqref{matrix-Y} satisfy normalization \eqref{asymptotics-Y}. 

Conversely, suppose that the matrix $Y(z)$ solves the $\delb$-problem \eqref{d-bar-Y} with the asymptotics 
\eqref{asymptotics-Y}. It follows from the special form \eqref{matrix-G} of the matrix $G(z)$ that $(Y_{n})_{11}(z)=P_{n}(z)$ and $(Y_{n})_{21}(z)=Q_{n-1}(z)$ --- polynomials of orders $n$ and $n-1$ respectively, 
and
\begin{align*}
\frac{\del}{\del\z}(Y_{n})_{12}(z) &=-\overline{P_{n}(z)}e^{-NV(z)}\chi_{D}(z), \\
\frac{\del}{\del\z}(Y_{n})_{22}(z) &=-\overline{Q_{n-1}(z)}e^{-NV(z)}\chi_{D}(z).
\end{align*}
Now it follows from normalization \eqref{asymptotics-Y} that
the leading coefficient of polynomial $P_{n}(z)$ is $1$ and
\begin{align*}
(Y_{n})_{12}(z) &= \frac{1}{\pi}\iint\limits_{D}\frac{\overline{P_{n}(z')}}{z'-z}e^{-NV(z')}d^{2}z', \\
(Y_{n})_{22}(z) & =\frac{1}{\pi}\iint\limits_{D}\frac{\overline{Q_{n-1}(z')}}{z'-z}e^{-NV(z')}d^{2}z'.
\end{align*}
Using geometric series and normalization \eqref{asymptotics-Y} once again, we obtain that polynomials $P_{n}$ and $Q_{n-1}$ satisfy
\begin{align*}
\iint\limits_{D}P_{n}(z)\z^{m}e^{-NV(z)}d^{2}z & =0 \quad\quad\text{for}\quad m<n,\\
-\iint\limits_{D}Q_{n-1}(z)\z^{m-1}e^{-NV(z)}d^{2}z & =\pi\delta_{mn}\quad\text{for}\quad m\leq n.
\end{align*}
From here it follows that 
$$\iint\limits_{D}P_{n}(z)\overline{P_{m}(z)} e^{-NV(z)}d^{2}z=0\quad\text{for all}\quad m<n,$$
which is sufficient to conclude that $P_{n}(z)$ are orthogonal polynomials on $\CC$ with the weight $e^{-NV(z)}\chi_{D}(z)$. Finally, polynomials $Q_{n}(z)$ satisfy
$$\iint\limits_{D}Q_{n-1}(z)\overline{P_{m-1}(z)}e^{-NV(z)}d^{2}z=-\pi\delta_{mn}\quad\text{for all}\quad m\leq n,$$
so that $\displaystyle{Q_{n}(z)=-\frac{\pi}{h_{n}}P_{n}(z)}$.
 \end{proof}
Similar to the case of the usual orthogonal polynomials, Proposition \ref{prop1}
reduces the asymptotic analysis of the orthogonal polynomials
(\ref{orthogonality}) to the asymptotic analysis of the  solution of
the $\delb$-problem (\ref{d-bar-Y})-(\ref{asymptotics-Y}).

\section{Elbau-Felder potential. Towards a normal matrix version of the DKMVZ asymptotic approach.} 
We will consider the matrix model with the weight $e^{-NV(z)}\chi_{D}(z)$, where $V(z)$ is the Elbau-Felder \cite{EF} potential  
\begin{equation} \label{E-F-potential}
V(z)=\frac{1}{t_{0}}\left(|z|^{2}-2\re\sum_{k=1}^{n+1}t_{k}z^{k}\right),
\end{equation}
where $t_{1}=0$, $|t_{2}|<1/2$ and $t_{0}V(z)$ is positive on $D\setminus\{0\}$.
Using again the Hermitian matrix model analogy, we shall expect that a fundamental role in the asymptotic analysis of the  $\delb$-problem (\ref{d-bar-Y})-(\ref{asymptotics-Y})
will be played by the {\it equilibrium measure}.
\begin{definition} An equilibrium measure for $V$ on $D$ is a Borel
probability measure $\mu$ on $D$ without point masses so that
$$
I(\mu) = \mbox{inf}\, I(\nu), \quad \nu \subset {\Bbb M}(D),
$$
where ${\Bbb M}(D)$ is the set of all  Borel
probability measures $\mu$ on $D$ without point masses, and the functional
$I(\nu)$ is defined by the equation,
$$
I(\nu):= \int V(z)d\nu(z) + \int \int _{z\neq \zeta} \log|z-\zeta|^{-1}d\nu(z)d\nu(\zeta).
$$
\end{definition}
\begin{theorem}[Elbau-Felder \cite{EF}]  \label{th7} There is $\delta>0$ such that for all
$0<t_{0}<\delta$ the unique equilibrium measure $d\mu$ exists and is given
by
$$d\mu(z)=\frac{1}{\pi t_{0}}\chi_{D_{+}}(z)d^{2}z,$$
where the domain $D_{+}\subset D$ contains the origin and has the property that
\begin{align*}
t_{0} & =\frac{1}{\pi}\iint\limits_{D_{+}}d^{2}z, \\
t_{k} & =-\frac{1}{\pi k}\iint\limits_{\CC\setminus D_{+}}z^{-k}d^{2}z=\frac{1}{2\pi ik}\oint\limits_{\del D_{+}}\z z^{-k}dz,\quad k=1,\dots,n+1,\\
t_{k} &=0,\quad j>n+1.
\end{align*}
These relations determines $D_{+}$ uniquely. In fact, the boundary $\Gamma$
of $D_{+}$ is a polynomial curve of degree $n$, i.e. $\Gamma$ is a smooth
simple closed curve in the complex plane with a parametrization $h:
S^{1}\subset {\Bbb C} \to {\Bbb C}$ of the form,
$$
h(w) = rw + a_{0} + a_{1}w^{-1} + ... + a_{n}w^{-n}, \quad |w| =1,
$$
with $r >0$ and $a_{n}\neq 0$. 
The equilibrium measure has the following properties. Set
$$E(z)=V(z)+2\iint\log|z-\zeta|^{-1}d\mu(\zeta).$$
\begin{itemize}
\item[\bf{1.}] $E(z)=E_{0}$ --- a constant --- for $z\in D_{+}$.
\item[\bf{2.}] $E(z) \geq E_{0}$ --- for $z \in D\setminus D_{+}$.
\end{itemize}
\end{theorem}
\subsection{A naive DKMVZ scheme.}
Suppose that there is an analytic function
$g(z)$ with the following properties.
\begin{itemize}
\item[\bf{1.}] $g(z)=\log z+O\left(\frac{1}{z}\right)$ as $|z|\rightarrow\infty$.
\item[\bf{2.}] $V(z)-g(z)-\overline{g(z)}=E_{0}$ on $D_{+}$.
\item[\bf{3.}] $V(z)-g(z)-\overline{g(z)}>E_{0}$ on $\CC\setminus D_{+}$.
\end{itemize}
Such function $g(z)$ could be used to study the asymptotics of the matrix $Y_{n}(z)$ in the limit,
\begin{equation}\label{limit}
n, N\rightarrow \infty, \quad \frac{n}{N}=\gamma, \quad \mbox{$\gamma$ is fixed},
\end{equation}
in exactly the same manner is it is done in the case of the
orthogonal polynomials in the line (see \cite{DKMVZ} and \cite{D2}). Namely, set
$$V_{\gamma}(z)=\frac{1}{\gamma}V(z)$$
and consider corresponding equilibrium measure $d\mu_{\gamma}$
(assuming that $0<\gamma t_{0}<\delta$) with the domain $D_{+}(\gamma)$ and the  corresponding function $g_{\gamma}(z)$ 
satisfying properties \textbf{1--3} (with $D_{+}$ and $E_{0}$ replaced by $D_{+}(\gamma)$ and $E_{0}(\gamma)$).
Then we can ``undress'' the $\delb$-problem  \eqref{d-bar-Y}--\eqref{matrix-G} with the normalization \eqref{asymptotics-Y} by setting
 $$Y_{n}(z)=e^{-\frac{nE_{0}(\gamma)}{2}\sigma_{3}}\Psi_{n}(z)e^{ng_{\gamma}(z)\sigma_{3} +\frac{nE_{0}(\gamma)}{2}\sigma_{3}}.$$
 The resulting matirx $\Psi_{n}(z)$ satisfies the simplified $\delb$-problem
 \begin{equation} \label{d-bar-Psi}
\frac{\del}{\del\z}\Psi_{n}(z)=\overline{\Psi_{n}(z)}
\begin{pmatrix} 0 & -e^{-n(V_{\gamma}(z)-g_{\gamma}(z)-\overline{g_{\gamma}(z)}-E_{0}(\gamma))} \\
0 & 0
\end{pmatrix}
 \end{equation}
with the standard normalization
\begin{equation} \label{asymptotics-Psi}
\Psi_{n}(z)=I +O\left(\frac{1}{z}\right)\quad\text{as}\quad |z|\rightarrow\infty,
\end{equation}
which follows from property \textbf{1} of the function $g_{\gamma}(z)$.
 
It is easy to pass to limit  (\ref{limit}) in the $\delb$-problem \eqref{d-bar-Psi}--\eqref{asymptotics-Psi}. Indeed, it follows from properties
\textbf{2-3} of the function $g_{\gamma}(z)$ that 
$$\lim_{n,N\rightarrow\infty}\Psi_{n}(z)=\Psi_{0}(z),$$
where the matrix $\Psi_{0}(z)$ satisfies the following \emph{model $\delb$-problem}.
\begin{equation} \label{d-bar-model}
\frac{\del}{\del\z}\Psi_{0}(z)=\overline{\Psi_{0}(z)}\begin{cases}
(\begin{smallmatrix} 0 & -1\\
0 & 0\end{smallmatrix}) & z\in D_{+}(\gamma) \\
0 & z\notin D_{+}(\gamma)
\end{cases}
\end{equation}
with the standard normalization
\begin{equation} \label{asymptotics-Psi-0}
\Psi_{0}(z)=I +O\left(\frac{1}{z}\right)\quad\text{as}\quad |z|\rightarrow\infty.
\end{equation}
This model $\delb$-problem is easily solved explicitly,
$$\Psi_{0}(z)=\begin{pmatrix} 1 & \frac{1}{\pi}\iint\limits_{D_{+}(\gamma)}\frac{1}{\zeta-z}d^{2}\zeta\\
0 & 1\end{pmatrix}.$$
\subsection{Function $g(z)$} 
Of course, the main assumption that there is an analytic function function $g(z)$ satisfying the properties \textbf{1-3} is not correct. Firstly, it follows from the property \textbf{1} that $g(z)$ in the neighborhood of infinity is defined up to an integer multiple of $2\pi i$, which is not a drawback since $e^{ng(z)\sigma_{3}}$ si well-defined for integer $n$. Secondly, the property \textbf{2} implies that the function
$V(z)$ is harmonic in $D_{+}$, which clearly contradicts \eqref{E-F-potential}. Adding to the confusion is the formal manipulation
$$\log|z-\zeta|^{2}= \log(z-\zeta) +\log(\z-\bar{\zeta}),$$
which suggests that
\begin{equation} \label{naive}
g(z)=\iint\limits\log(z-\zeta)d\mu(\zeta)=\frac{1}{\pi t_{0}}\iint\limits_{D_{+}}\log(z-\zeta)d^{2}\zeta
\end{equation} 
satisfies properties \textbf{1-3}. However, this is not so since we need to treat carefully the branches of $\log$ in order to define the integral in \eqref{naive} and investigate its analytic properties.

For this aim, consider the logarithmic potential given by the uniform distribution of charges in the domain $D$,
$$V_{0}(z)=\iint\limits_{D}\log|z-w|^{2}d^{2}w.$$
Let $\Gamma=\del D$ with fixed point $\zeta_{0}\in\Gamma$. 
For $z\in D$ denote by $D_{\vep}(z)$ domain obtained by removing the disk of radius $\vep$ around $z$, so that 
$$\del D_{\vep}(z)=\Gamma\cup-C_{\vep}(z),$$
where $C_{\vep}(z)$ is the circle $|w-z|=\vep$ oriented counter-closkwise, and the minus sign denotes negative orientation. 
Since
$$\log|w-z|^{2}dw\wedge d\w=-d(\log|w-z|^{2}\w dw)-\frac{\w}{\w-\z}\,dw\wedge d\w,$$
by Stokes' theorem, we have
\begin{align*}
V_{0}(z) & =\frac{i}{2}\iint\limits_{D}\log|z-w|^{2}dw\wedge d\bar{w} \\
&=\frac{i}{2}\lim_{\vep\rightarrow 0}\iint\limits_{D_{\vep}(z)}\left(
-d(\log|w-z|^{2}\w dw)-\frac{\w}{\w-\z}\,dw\wedge d\w \right) \\
&=\frac{1}{2i}\lim_{\vep\rightarrow 0}\iint\limits_{D_{\vep}(z)}
d\left(\log|w-z|^{2}\w dw+\frac{\w}{\w-\z}\,wd\w\right) \\
&=\frac{1}{2i}\lim_{\vep\rightarrow 0}\oint\limits_{\del D_{\vep}(z)}
\left(\log|\zeta-z|^{2}\bar\zeta d\zeta+\frac{\bar\zeta}{\bar\zeta-\z}\,\zeta d\bar\zeta\right).
\end{align*}
Now 
\begin{align*}
\lim_{\vep\rightarrow 0}\int\limits_{C_{\vep}(z)}
\log|\zeta-z|^{2}\bar\zeta d\zeta = 0,\quad\quad
\lim_{\vep\rightarrow 0}\int\limits_{C_{\vep}(z)}\frac{\bar\zeta}{\bar\zeta-\z}\,\zeta d\bar\zeta = -2\pi i|z|^{2},
\end{align*}
so that
\begin{align*}
V_{0}(z) & =\pi |z|^{2} + \frac{1}{2i}\oint\limits_{\Gamma}
\left(\log|\zeta-z|^{2}\bar\zeta d\zeta+\frac{\bar\zeta}{\bar\zeta-\z}\,\zeta d\bar\zeta\right).
\end{align*}
Set $\omega=\w dw$ and define a function  
$\Omega$ on $\Gamma\setminus\{\zeta_{0}\}$ by 
$$\Omega(\zeta)=\int\limits_{\zeta_{0}}^{\zeta}\omega,$$
where the integration is along the oriented path in $\Gamma$ connecting points $\zeta_{0}$ and $\zeta$. We have $\Omega_{-}(\zeta_{0})=0$ for a path consisting of a single point $\zeta_{0}$, and
$$\Omega_{+}(\zeta_{0})=\oint\limits_{\Gamma}\w dw=\iint\limits_{D}d\w\wedge dw=2i A(D),$$
for the loop $\Gamma$ starting and ending at $\zeta_{0}$,
where $A(D)$ is the area of $D$. Thus
\begin{align*}
\oint\limits_{\Gamma} \log|\zeta-z|^{2}\zeta d\zeta & =\oint\limits_{\Gamma}
\log|\zeta-z|^{2}d\Omega(\zeta) =\left.\Delta(\log|z-\zeta|^{2}\Omega(\zeta))\right|_{\zeta_{0}}^{\zeta_{0}} \\
&-\oint\limits_{\Gamma}\Omega(\zeta)\left(\frac{d\zeta}{\zeta-z} +
\frac{d\bar\zeta}{\bar\zeta-\z}\right) \\
& =2iA(D)\log|z-\zeta_{0}|^{2}-\oint\limits_{\Gamma}\Omega(\zeta)\left(\frac{d\zeta}{\zeta-z} +
\frac{d\bar\zeta}{\bar\zeta-\z}\right),
\end{align*}
so that
$$V_{0}(z)=\pi|z|^{2} +A(D)\log|z-\zeta_{0}|^{2} +\frac{i}{2}\oint\limits_{\Gamma}\left(\frac{\Omega(\zeta)}{\zeta-z}\,d\zeta +
\frac{\Omega(\zeta)}{\bar\zeta-\z}\,d\bar\zeta -\frac{\bar\zeta}{\bar\zeta-\z}\,\zeta d\bar\zeta \right).$$
Since the potential $V_{0}$ is real-valued, we have
\begin{align*}
V_{0}(z) & =\pi|z|^{2} +A(D)\log|z-\zeta_{0}|^{2} +\frac{i}{4}\oint\limits_{\Gamma}\left(\frac{\Omega(\zeta)}{\zeta-z}\,d\zeta -\frac{\overline{\Omega(\zeta)}}{\bar\zeta-\z}\,d\bar\zeta \right.\\
 & \left.-\frac{\overline{\Omega(\zeta)}}{\zeta-z}\,d\zeta +
\frac{\Omega(\zeta)}{\bar\zeta-\z}\,d\bar\zeta +\frac{\zeta}{\zeta-z}\,\bar\zeta d\zeta -\frac{\bar\zeta}{\bar\zeta-\z}\,\zeta d\bar\zeta \right) \\
&=\pi|z|^{2} +A(D)\log|z-\zeta_{0}|^{2} +\frac{i}{4}\oint\limits_{\Gamma}\left(\frac{\Omega(\zeta)-\overline{\Omega(\zeta)} +|\zeta|^{2} }{\zeta-z}\,d\zeta  \right. \\
& -\left. \frac{\overline{\Omega(\zeta)} - \Omega(\zeta) +|\zeta|^{2} }{\bar\zeta-\z}\,d\bar\zeta\right). 
\end{align*}
Finally, observing that $\omega +\bar\omega =d|w|^{2}$, we get
$$\Omega(\zeta)+\overline{\Omega(\zeta)}=\int\limits_{\zeta_{0}}^{\zeta}d|w|^{2}=|\zeta|^{2}-|\zeta_{0}|^{2},$$
so that 
$$\Omega(\zeta)-\overline{\Omega(\zeta)} +|\zeta|^{2}=2\Omega(\zeta) +|\zeta_{0}|^{2},$$
and we obtain 
\begin{equation} \label{g-inside}
V_{0}(z)=\pi(|z|^{2}-|\zeta_{0}|^{2}) +A(D)\log|z-\zeta_{0}|^{2} +
\frac{i}{2}\oint\limits_{\Gamma}\left(\frac{\Omega(\zeta)}{\zeta-z}\,d\zeta - \frac{\overline{\Omega(\zeta)}}{\bar\zeta-\z}\,d\bar\zeta\right),
\end{equation}
where $z\in D$. This is a desired representation of the area potential $V_{0}(z)$ as the real part of the first derivative of a single layer potential. 

The same computation for $z\in\CC\setminus\bar{D}$ gives
\begin{equation} \label{g-outside}
V_{0}(z)=A(D)\log|z-\zeta_{0}|^{2} +
\frac{i}{2}\oint\limits_{\Gamma}\left(\frac{\Omega(\zeta)}{\zeta-z}\,d\zeta - \frac{\overline{\Omega(\zeta)}}{\bar\zeta-\z}\,d\bar\zeta\right),
\end{equation}

Returning to Elbau-Felder potential $V(z)$ and setting $D=D_{+}$,
$\Gamma=\del D_{+}$,
we get $\displaystyle{E(z) = V(z) -\frac{1}{\pi t_{0}}V_{0}(z)}$, so that 
\begin{equation}\label{constant}
V(z) -\frac{1}{t_{0}}(|z|^{2}-|\zeta_{0}|^{2}) - \log|z-\zeta_{0}|^{2} -\frac{i}{2\pi t_{0}}\oint\limits_{\Gamma}\left(\frac{\Omega(\zeta)}{\zeta-z}\,d\zeta - \frac{\overline{\Omega(\zeta)}}{\bar\zeta-\z}\,d\bar\zeta\right)=E_{0}
\end{equation}
when $z\in D_{+}$, and
\begin{equation}\label{nonconstant}
V(z) -\log|z-\zeta_{0}|^{2} -\frac{i}{2\pi t_{0}}\oint\limits_{\Gamma}\left(\frac{\Omega(\zeta)}{\zeta-z}\,d\zeta - \frac{\overline{\Omega(\zeta)}}{\bar\zeta-\z}\,d\bar\zeta\right)=E(z)
\end{equation}
 when $z\in D_{-}$.
\begin{remark} We note  that equation (\ref{constant}), i.e. the statement that
the l.h.s. of (\ref{constant}) is constant when $z\in D_{+}$,  is equivalent to
the moment equations of Theorem \ref{th7} which determine the contour
$\Gamma$ (cf.\cite{EF}, p.12, Lemma 6.3).
\end{remark}

  Now we are ready to introduce the function $g(z)$. Namely, set
\begin{equation} \label{g-definition}
g(z)=\log(z-\zeta_{0}) + \frac{i}{2\pi t_{0}}\oint\limits_{\Gamma}\frac{\Omega(\zeta)}{\zeta-z}\,d\zeta.
\end{equation}
The function $g(z)$ is holomorphic in $\CC\setminus\Gamma$, is multi-valued with periods $2\pi i\mathbb{Z}$ 
(singe-valued on the plane with the outside cut starting from $\zeta_{0}$) and has the asymptotics
$$g(z)=\log z+O(z^{-1})\quad\text{as}\quad z\rightarrow\infty.$$
The function $e^{ng(z)}$ is single-valued for $n\in\mathbb{Z}$.
The function $g(z)$ is discontinuous on $\Gamma$ (by Sokhotski-Plemelj formula).

We summarize this as the following statement.
\begin{proposition} The Elbau-Felder potential $V(z)$ has the following
representations
\begin{itemize}
\item[(i)] For $z\in D_{+}$,
$$V(z) - g(z) -\overline{g(z)} =E_{0} + \frac{1}{t_{0}}(|z|^{2}-|\zeta_{0}|^{2}).$$
\item[(ii)] For $z\in D_{-}$,
$$V(z) - g(z) -\overline{g(z)} =E(z).$$
\item[(iii)] For $z\in D\setminus D_{+}$,
\begin{equation}\label{ineq}
V(z) - g(z) -\overline{g(z)} =E(z) > E_{0}.
\end{equation}
\end{itemize}
\end{proposition}

\subsection{A first possible version of the DKMVZ scheme}
The correct strategy is now the following. Let $g_{\gamma}(z)$,
$D_{+}(\gamma)$, $E_{0}(\gamma)$, etc. denote the respective
objects associated with the potential $V_{\gamma}(z)$. We set
\begin{equation}\label{Psi-cor}
 Y_{n}(z)=e^{-\frac{nE_{0}(\gamma)}{2}\sigma_{3}+\frac{n|\zeta_{0}|^2}{2\gamma t_{0}}\sigma_{3}}
\Psi_{n}(z)e^{ng_{\gamma}(z)\sigma_{3} +\frac{nE_{0}(\gamma)}{2}\sigma_{3}
-\frac{n|\zeta_{0}|^2}{2\gamma t_{0}}\sigma_{3}}.
\end{equation}
 The resulting matirx $\Psi_{n}(z)$ satisfies the  $\delb$-problem (the correct version of (\ref{d-bar-Psi}))
 $$
 \frac{\del}{\del\z}\Psi_{n}(z)=\overline{\Psi_{n}(z)}
\begin{pmatrix} 0 & -e^{- \frac{n|z|^2}{\gamma t_{0}}} \\
0 & 0
\end{pmatrix}, \quad z \in D_{+},
$$
 \begin{equation} \label{d-bar-Psi-cor}
\frac{\del}{\del\z}\Psi_{n}(z)=\overline{\Psi_{n}(z)}
\begin{pmatrix} 0 & -e^{-n\left(E(z)-E_{0}(\gamma)+ \frac{|\zeta_{0}|^2}{\gamma t_{0}}\right)}\chi_{D_{+}}(z) \\
0 & 0
\end{pmatrix}, \quad z\in \CC\setminus D_{+},
 \end{equation}

$$
\Psi_{n+}(z) = \Psi_{n-}(z)e^{\frac{n}{\gamma t_{0}}\Omega(z)\sigma_{3}}, \quad z \in
 \Gamma \equiv \Gamma(\gamma),
$$
with the standard normalization
\begin{equation} \label{asymptotics-Psi-cor}
\Psi_{n}(z)=I +O\left(\frac{1}{z}\right)\quad\text{as}\quad |z|\rightarrow\infty.
\end{equation}
By virtue of  condition (\ref{ineq}), we expect 
that the limiting function $\Psi^{0}_{n}(z)$
satisfies the model problem
\begin{equation} \label{d-bar-model-cor}
\frac{\del}{\del\z}\Psi^{0}_{n}(z)=\overline{\Psi^{0}_{n}(z)}\begin{cases}
\begin{pmatrix} 0 & -e^{- \frac{n|z|^2}{\gamma t_{0}}}\\
0 & 0\end{pmatrix} & z\in D_{+}(\gamma) \\\\
0 & z\notin D_{+}(\gamma)
\end{cases}
\end{equation}

$$
\Psi^{0}_{n+}(z) = \Psi^{0}_{n-}(z)e^{\frac{n}{\gamma t_{0}}\Omega(z)\sigma_{3}}, \quad z \in \Gamma,
$$
with the standard normalization
\begin{equation} \label{asymptotics-Psi-0-cor}
\Psi_{0}(z)=I +O\left(\frac{1}{z}\right)\quad\text{as}\quad |z|\rightarrow\infty.
\end{equation}
The open questions are now the following.
\begin{enumerate}
\item How to solve this model problem?
The ``unfortunate'' thing is the presence of complex
conjugation in (\ref{d-bar-model-cor}). Indeed, if we neglect the jump across the
contour $\Gamma$, the 
$\delb$-problem alone can be of course solved explicitly,
\begin{equation}\label{final}
\Psi_{0}(z)=\begin{pmatrix} 1 &\frac{1}{\pi}\iint\limits_{D_{+}(\gamma)}
\frac{e^{-\frac{n|\zeta|^2}{\gamma t_{0}}}}{\zeta-z}d^{2}\zeta\\
0 & 1\end{pmatrix},
\end{equation}
and the solution won't have any jumps.  If not for the complex conjugation,
the function $\Psi_{0}(z)$ could be used to undress in the usual way
problem (\ref{d-bar-model-cor}) and reduce it to a pure Riemann-Hilbert problem.
\item The arguments that led us to the model problem (\ref{d-bar-model-cor})
and which are based on inequality (\ref{ineq}),
even on the formal level, are not very convincing: the real part of $\Omega(z)$ is $|z|^2 - |\zeta_{0}|^2\neq 0$ so that 
the diagonal jump matrix on $\Gamma$ is not pure oscillatory.
\end{enumerate}

\subsection{A second possible version of the DKMVZ scheme}
The above deficiency of the proposed analog of the DKMVZ scheme can be partially overcome by performing the following modification. Let us replace the
function $\Omega(\zeta)$ by the function,
$$
\Omega_{0}(\zeta) =\frac{1}{2} \int_{\zeta_{0}}^{\zeta}(\w dw -wd\w).
$$
The function $\Omega_{0}(\zeta)$ is pure imaginary on $\Gamma$ and
it is related to the function $\Omega(\zeta)$ by the equation,
$$
\Omega_{0}(\zeta) = \Omega(\zeta) -\frac{|\zeta|^{2}  - |\zeta_{0}|^{2}}{2}.
$$
By using again Stokes' theorem, we observe that
$$
\frac{i}{4\pi t_{0}}\oint\limits_{\Gamma}(|\zeta|^2 - |\zeta_{0}|^2)
\left(\frac{d\zeta}{\zeta-z} - \frac{d\bar\zeta}{\bar\zeta-\z}\right)
$$
$$
= - \frac{1}{2\pi t_{0}}\iint_{D_{+}}\left(\frac{\zeta}{\zeta - z}
+\frac{\bar\zeta}{\bar\zeta - \z}\right)d^2\zeta -\frac{1}{t_{0}}
(|z|^2 - |\zeta_{0}|^2).
$$
This allows to re-write equations (\ref{constant}) and (\ref{nonconstant})
in the form,
$$
V(z)  - \log|z-\zeta_{0}|^{2} -\frac{i}{2\pi t_{0}}\oint\limits_{\Gamma}\left(\frac{\Omega_{0}(\zeta)}{\zeta-z}\,d\zeta - \frac{\overline{\Omega_{0}(\zeta)}}{\bar\zeta-\z}\,d\bar\zeta\right)
$$
\begin{equation}\label{constant1}
+\frac{1}{2\pi t_{0}}\iint_{D_{+}}\left(\frac{\zeta}{\zeta - z}
+\frac{\bar\zeta}{\bar\zeta - \z}\right)d^2\zeta=E_{0}
\end{equation}
when $z\in D_{+}$, and
$$
V(z)  - \log|z-\zeta_{0}|^{2} -\frac{i}{2\pi t_{0}}\oint\limits_{\Gamma}\left(\frac{\Omega_{0}(\zeta)}{\zeta-z}\,d\zeta - \frac{\overline{\Omega_{0}(\zeta)}}{\bar\zeta-\z}\,d\bar\zeta\right)
$$
\begin{equation}\label{constant2}
+\frac{1}{2\pi t_{0}}\iint_{D_{+}}\left(\frac{\zeta}{\zeta - z}
+\frac{\bar\zeta}{\bar\zeta - \z}\right)d^2\zeta=E(z)
\end{equation}
 when $z\in D_{-}$.
These formulae in turn yield the following modification
of the definition (\ref{g-definition}) of the $g$-function
\begin{equation} \label{g-definition1}
g(z)=\log(z-\zeta_{0}) + \frac{i}{2\pi t_{0}}\oint\limits_{\Gamma}\frac{\Omega_{0}(\zeta)}{\zeta-z}\,d\zeta
- \frac{1}{2\pi t_{0}}\iint_{D_{+}}\frac{\zeta}{\zeta - z}d^2\zeta.
\end{equation}
Note that the function $g(z)$ is not holomorphic in ${\Bbb C}\setminus {\Gamma}$ anymore! In fact\footnote{Probably this property of the function $g(z)$ reflects a major difference between the Riemann-Hilbert problem and the $\delb$-problem.},
\begin{equation}\label{gextra}
\frac{\del}{\del\z}g(z) = \frac{z}{2 t_{0}}\chi_{D_{+}}(z).
\end{equation}
A slight modification is also needed in the definition of the function $\Psi_{n}(z)$;
indeed, we should put,
\begin{equation}\label{Psi-cor1}
 Y_{n}(z)=e^{-\frac{nE_{0}(\gamma)}{2}\sigma_{3}}
\Psi_{n}(z)e^{ng_{\gamma}(z)\sigma_{3} +\frac{nE_{0}(\gamma)}{2}\sigma_{3}
}.
\end{equation}
Taking into account (\ref{gextra}), the $\delb$-problem for the matrix $\Psi$ now reads,
 \begin{equation}\label{dbarPsi2}
 \frac{\del}{\del\z}\Psi_{n}(z) + n\frac{z}{2 \gamma t_{0}}\Psi_{n}(z)\sigma_{3}
 =\overline{\Psi_{n}(z)}
\begin{pmatrix} 0 & - 1 \\
0 & 0
\end{pmatrix}, \quad z \in D_{+},
\end{equation}
 \begin{equation} \label{d-bar-Psi-cor1}
\frac{\del}{\del\z}\Psi_{n}(z)=\overline{\Psi_{n}(z)}
\begin{pmatrix} 0 & -e^{-n\left(E(z)-E_{0}(\gamma)\right)}\chi_{D_{+}}(z) \\
0 & 0
\end{pmatrix}, \quad z\in \CC\setminus D_{+},
 \end{equation}

$$
\Psi_{n+}(z) = \Psi_{n-}(z)e^{\frac{n}{\gamma t_{0}}\Omega_{0}(z)\sigma_{3}}, \quad z \in
 \Gamma \equiv \Gamma(\gamma),
$$
with the standard normalization
\begin{equation} \label{asymptotics-Psi-cor}
\Psi_{n}(z)=I +O\left(\frac{1}{z}\right)\quad\text{as}\quad |z|\rightarrow\infty.
\end{equation}

The function $\Omega_{0}(z)$ is now purely imaginary. Therefore,
the arguments based on inequality (\ref{ineq}) seem to be
more sound than in the previous approach, and they lead us to
the following new model $\delb$-problem
\begin{equation}\label{newmod}
 \frac{\del}{\del\z}\Psi^{0}_{n}(z) + n\frac{z}{2 \gamma t_{0}}\Psi^{0}_{n}(z)\sigma_{3}
 =\overline{\Psi^{0}_{n}(z)}
\begin{pmatrix} 0 & - 1 \\
0 & 0
\end{pmatrix}, \quad z \in D_{+},
\end{equation}
 \begin{equation} \label{newmod2}
\frac{\del}{\del\z}\Psi^{0}_{n}(z)=0, \quad z\in \CC\setminus D_{+},
 \end{equation}

$$
\Psi^{0}_{n+}(z) = \Psi^{0}_{n-}(z)e^{\frac{n}{\gamma t_{0}}\Omega_{0}(z)\sigma_{3}}, \quad z \in
 \Gamma \equiv \Gamma(\gamma),
$$
with the standard normalization
\begin{equation} \label{asymptotics-Psi-cor1}
\Psi^{0}_{n}(z)=I +O\left(\frac{1}{z}\right)\quad\text{as}\quad |z|\rightarrow\infty.
\end{equation}

\begin{remark} Due to the presence of the large parameter $n$ in the
left hand side of equation (\ref{dbarPsi2}), the transition to the
model problem (\ref{newmod})-(\ref{asymptotics-Psi-cor1}) is
still not quite satisfactory even on the formal level.
\end{remark}

\subsection{An important concluding remark}
In context of the theory of orthogonal
polynomials, a matrix $\delb$-problem has also appeared
in the recent work of K. McLaughlin and P. Miller \cite{MM} devoted
to orthogonal polynomials on the unite circle with the non-analytic
weights. However, unlike the problem (\ref{d-bar-Y})-(\ref{asymptotics-Y}), the
$\delb$-problem of \cite{MM} is not
the starting point of the analysis; indeed, the staring point of  \cite{MM}
is still the usual matrix Riemann-Hilbert problem and the $\delb$-problem
of  McLaughlin and  Miller  is introduced out of the necessity to modify
the ``openning lenses''  step of the usual DKMVZ scheme. Even more important
difference between the $\delb$-problem considered here and the $\delb$-problem in \cite{MM} is the absence of the complex conjugation in the 
basic $\delb$-relation. In  one hand, this fact simplifies the 
implementation of the ``undressing procedures'' --- the very important
technical element of all integrable asymptotic schemes. On the other
hand, as we have shown, the presence of the complex conjugation
in the right hand side of  \eqref{d-bar-Y} is truly essential for the
incorporation into the asymptotic analysis of the concepts of  
equilibrium measure and $g$-function.

In spite of these differences 
we believe that using methods of \cite{MM} will help allow
to overcome the indicated above obstacles in the asymptotic
analysis of the $\delb$-problem (\ref{d-bar-Y})-(\ref{asymptotics-Y}).
Specifically, we we think that one needs to develop and then apply
to the model problems  (\ref{d-bar-model-cor})-(\ref{asymptotics-Psi-0-cor}) or (\ref{newmod})-(\ref{asymptotics-Psi-cor1}) of the 
the $\delb$-version of the ``openning lenses''  step in the DKMVZ method. 

In conclusion, we want to point out at the $\delb$-method in the theory of integrable systems, introduced long ago by A.S. Fokas and M.J. Ablowitz in \cite{FA}, as yet another source of tools for the analysis of the $\delb$-problem \eqref{d-bar-Y}-\eqref{asymptotics-Y}. 
%

\section*{Acknowledgements}
Alexander Its was supported in part by the NSF grants DMS-0401009
and DMS-0701768, and Leon Takhtajan --- by the NSF grants DMS-0204628 and DMS-0705263.


\bigskip
\begin{thebibliography}{99}



\bibitem{D2} P.~A.~Deift, Orthogonal Polynomials and Random Matrices:
A Riemann-Hilbert Approach, {\it Courant Lecture Notes in Mathematics},
{\bf 3}, {\em CIMS, New York} (1999).



\bibitem {DKMVZ} 
P. Deift, T. Kriecherbauer, K. T-R. McLaughlin,
S. Venakides, and X. Zhou,
Uniform asymptotics for polynomials orthogonal with respect
to varying exponential weights and applications to universality
questions in random matrix theory,
{\it Commun. Pure Appl. Math.}, {\bf 52} (1999), 1335-1425. 

\bibitem{EF} P. Elbau and G. Felder, Density of eigenvalues of random
normal matrices, {\it Commun. Math. Phys}, {\bf 259} (2005), 433-450. 

\bibitem {FIK1} A.S. Fokas, A.R. Its and A.V. Kitaev, The isomonodromy approach
to matrix problems in 2D quantum gravity, {\it Commun.\ Math.\
Phys.}, {\bf 147} (1992), 395--430 .

\bibitem{FA} A.S. Fokas and M.J. Ablowitz, On the inverse scattering transform of multidimensional nonlinear equations related to first-order systems in the plane, {\it J. Math. Phys.}, {\bf 25} (1984), 2494-2505.

\bibitem{MM} K. T.-R. McLaughlin and P. D. Miller,
 The $\delb$ steepest descent method and the asymptotic behavior
 of polynomials orthogonal on the unit circle with
 fixed and exponentially varying nonanalytic weights,
 {\it IMPR}, v. 2006 (2006), 1-78.
\bibitem {WZ} P. B. Wiegmann, A. Zabrodin, Conformal maps and integrable
hierarchise, {\it Commun. Math. Phys.},  {\bf 213} (2000), 523-538.
\bibitem{KKMWWZ} I. K. Kostov, I. Krichever, M. Mineev-Weinstein, P. B. Wiegmann, A. Zabrodin,  $\tau$-function for analytic curves, {\it Random matrix models and their applications},
285-299, Math. Sci. Res. Inst. Publ., {\bf 40}, Cambridge Univ. Press, Cambridge, 2001.
\end{thebibliography}
 \end{document}